\documentclass{article}

\usepackage{amsthm}
\usepackage{amssymb}
\usepackage{amsxtra}

\newtheorem{thm}{Theorem}[section]

\newtheorem{conj}{Conjecture}[section]

\title{Higher moments of primes in short intervals II}
\author{Tsz Ho Chan}

\begin{document}
\maketitle
%-----------------------------------------------------------------
\begin{abstract}
Given good knowledge on the even moments, we derive asymptotic
formulas for $\lambda$-th moments of primes in short intervals and
prove ``equivalence'' result on odd moments. We also provide
numerical evidence in support of these results.
\end{abstract}
%-----------------------------------------------------------------
\section{Introduction}

In [\ref{MS}], Montgomery and Soundararajan studied the moments
$$M_k(X;h) := \int_{1}^{X} (\psi(x+h) - \psi(x) - h)^k dx$$
where $k$ is a positive integer, $\psi(x) = \sum_{n \leq x}
\Lambda(n)$ and $\Lambda(n)$ is von Mangoldt lambda function. They
proved that, under a strong form of Hardy-Littlewood prime-$k$
tuple conjecture, for small $\epsilon > 0$, there is a $\epsilon'
> 0$,
\begin{equation}
\label{1.1} M_k(X;h) = \mu_k h^{k/2} \int_{1}^{X}
(\log{\frac{x}{h}} + B)^{k/2} dx + O_k(h^{k/2} X^{1-\epsilon})
\end{equation}
uniformly for $(\log X)^{15 k^2} \leq h \leq X^{1/k - \epsilon'}$
where $\mu_k = 1 \cdot 3 \cdot \cdot \cdot (k-1)$ if $k$ is even,
and $\mu_k = 0$ if $k$ is odd. Here $B = 1 - C_0 - \log 2\pi$ and
$C_0$ denotes Euler's constant. One further expects that
(\ref{1.1}) holds uniformly for $X^{\epsilon'} \leq h \leq
X^{1-\epsilon'}$. Consider the closely related moments:
\begin{equation}
\label{1.2} \widetilde{M}_k(X;\delta) := \int_{1}^{X} (\psi(x +
\delta x) - \psi(x) - \delta x)^k dx.
\end{equation}
In [\ref{C}], the author proved that, roughly speaking,
(\ref{1.1}) holding uniformly for $X^{\epsilon'} \leq h \leq
X^{1-\epsilon'}$ is equivalent to
\begin{equation}
\label{1.3} \widetilde{M}_k(X;\delta) =
\frac{\mu_k}{\frac{k}{2}+1} X^{k/2+1} \delta^{k/2}
\Bigl(\log{\frac{1}{\delta}} + B \Bigr)^{k/2} + O_k(\delta^{k/2}
X^{k/2 + 1 - \epsilon})
\end{equation}
holding uniformly for $X^{-1+\epsilon'} \leq \delta \leq
X^{-\epsilon'}$.

In this article, we shall study the following more general
moments: For $\lambda > 0$,
\begin{equation}
\label{1.4} M_\lambda(X;h) := \int_{1}^{X} |\psi(x+h) - \psi(x) -
h|^\lambda dx,
\end{equation}
\begin{equation}
\label{1.5} \widetilde{M}_\lambda(X;\delta) := \int_{1}^{X}
|\psi(x + \delta x) - \psi(x) - \delta x|^\lambda dx.
\end{equation}
Instead of (\ref{1.1}) and (\ref{1.3}), we can assume the
following weaker versions: For some $\epsilon > 0$ and all even
positive integer $k$,
\begin{equation}
\label{1.6}
\begin{split}
\int_{1}^{X} (\psi(x+h) - \psi(x) - h)^k dx =&
\frac{\Gamma(k+1)}{\Gamma(\frac{k}{2}+1) 2^{k/2}} X h^{k/2}
\Bigl(\log{\frac{X}{h}}\Bigr)^{k/2} \\
&+ O(A^{B^k} X h^{k/2} \log^{k/2 - 1} X)
\end{split}
\end{equation}
for $X^\epsilon \leq h \leq X^{1-\epsilon}$, and
\begin{equation}
\label{1.7}
\begin{split}
\int_{1}^{X} (\psi(x+\delta x) - \psi(x) - \delta x)^k dx =&
\frac{\Gamma(k+1)}{\Gamma(\frac{k}{2}+2) 2^{k/2}} X^{k/2 + 1}
\delta^{k/2} \Bigl(\log{\frac{1}{\delta}}\Bigr)^{k/2} \\
&+ O(A^{B^k} X^{k/2+1} \delta^{k/2} \log^{k/2-1} X)
\end{split}
\end{equation}
for $X^{-1+\epsilon} \leq \delta \leq X^{-\epsilon}$. Remarks: 1.
$\Gamma(x)$ is the gamma function ($\Gamma(n+1) = n!$ for
non-negative integer $n$). 2. The weaker versions suffice because
our proof of the following theorems gives such poor error terms
that only the first main terms matter. 3. We make the implicit
constants' dependence on $k$ explicit in the error terms. Here
$A$, $B$ are some absolute constants greater than $1$ and
$A^{B^k}$ is the result of tracing the $k$ dependency of the error
terms in [\ref{MV}] and [\ref{MS}] explicitly. 4. $A$ and $B$ may
depend on $\epsilon$ but we treat $\epsilon$ as fixed to start
with.

We shall prove the following
\begin{thm}
\label{theorem1} If (\ref{1.6}) is true for every even positive
integer $k$ in $X^\epsilon \leq h \leq X^{1-\epsilon}$, then for
any $\lambda > 0$,
\begin{eqnarray*}
\int_{1}^{X} |\psi(x+h) - \psi(x) -h|^{\lambda} dx &=&
{\Gamma(\lambda+1) \over \Gamma({\lambda \over 2}+1)
2^{\lambda/2}} X h^{\lambda/2} (\log{X \over h})^ {\lambda/2} \\
&+& O_{\lambda} \Bigl({X h^{\lambda/2} (\log{X \over
h})^{{\lambda/2} } \over \log_{4}{X}} \Bigr)
\end{eqnarray*}
for $X^{\epsilon} \leq h \leq X^{1-\epsilon}$. Here $\log_{n}{x}$
stands for $n$ times iterated logarithm ($\log_{1}{x}=\log{x}$ and
$\log_{i+1}{x} = \log{\log_{i}{x}}$).
\end{thm}

\begin{thm}
\label{theorem2} If (\ref{1.7}) is true for every even positive
integer $k$ in $X^{-1+\epsilon} \leq \delta \leq X^{-\epsilon}$,
then for any $\lambda > 0$,
\begin{eqnarray*}
\int_{1}^{X} |\psi(x+\delta x) - \psi(x) - \delta x|^{\lambda} dx
&=& {\Gamma (\lambda+1) \over \Gamma({\lambda \over 2} + 2)
2^{\lambda/2}} X^{{ \lambda/2}+1} \delta^{\lambda/2} \Bigl(\log{1
\over \delta}\Bigr)^{\lambda/2}
\\
&+& O_{\lambda}\Bigl({X^{{\lambda/2}+1} \delta^{\lambda/2} (\log{1
\over \delta})^{\lambda/2} \over \log_{4}{X}} \Bigr)
\end{eqnarray*}
for $X^{-1+\epsilon} \leq \delta \leq X^{-\epsilon}$.
\end{thm}
Using these theorems, we shall prove ``equivalence'' result in
section 3. We shall give numerical evidence in the last section.
This work is part of the author's $2002$ PhD thesis.
%------------------------------------------------------------------
\section{Proof of Theorem \ref{theorem1} and \ref{theorem2}}

We need Legendre's double formula for gamma function:
\begin{equation}
\label{2.1} \sqrt{\pi} \Gamma(2z) = 2^{2z-1} \Gamma(z) \Gamma(z +
\frac{1}{2}).
\end{equation}
Since the proof of Theorem \ref{theorem1} and Theorem
\ref{theorem2} are very similar, we shall give the proof of
Theorem \ref{theorem1} only.

\bigskip

Proof of Theorem \ref{theorem1}: The method is essentially that of
Ghosh [\ref{G1}] and [\ref{G2}]. There are two cases depending on
the size of $\lambda$:

(i) $0 < \lambda \leq 1$.

Let $F(x) = {\psi(x+h)-\psi(x)-h \over 2^{-1/2} h^{1/2}
(\log{X/h})^{1/2}}$ and $G_{\lambda} = \int_{0}^{\infty}
{(\sin{u})^2 \over u^{1+\lambda}} du$. For any $\mu \geq 1$,
\begin{eqnarray*}
|F(x)|^{\lambda} &=& {1 \over G_{\lambda}} \int_{0}^{\infty}
{(\sin{|F(x)|u})^2
\over u^{1+\lambda} } du \\
&=&{1 \over G_{\lambda}} \int_{0}^{\mu} {(\sin{|F(x)|u})^2 \over
u^{1+\lambda}} du + O\Bigl({1 \over \lambda G_{\lambda}}
\mu^{-\lambda} \Bigr).
\end{eqnarray*}
Hence,
\begin{equation}
\label{2.2} \int_{1}^{X} |F(x)|^{\lambda} dx = {1 \over
G_{\lambda}} \int_{0}^{\mu} \int_{1} ^{X} {(\sin{|F(x)|u})^2 \over
u^{1+\lambda} } dx du + O_{\lambda}\Bigl({X \over
\mu^{\lambda}}\Bigr).
\end{equation}
Note that $(\sin{x})^2 = {1 \over 2}\sum_{j=1}^{N}
{(-1)^{j+1}(2x)^{2j} \over (2j)!} + O({(2x)^{2N+2} \over
(2N+2)!})$. So the main term of (\ref{2.2}) is
\begin{equation}
\label{2.3}
\begin{split}
& {1 \over 2G_{\lambda}} \int_{0}^{\mu} {1 \over u^{1+\lambda}}
\sum_{j=1}^{N} {(-1)^{j+1} (2u)^{2j} \over (2j)!}
\Bigl(\int_{1}^{X} |F(x)|^{2j} dx\Bigr) du
\\
+& O\Bigl({4^N \over (2N+2)!} \int_{0}^{\mu} u^{2N} du
\int_{1}^{X} |F(x)|^ {2N+2} dx \Bigr).
\end{split}
\end{equation}
Using assumption (\ref{1.6}), the error term in (\ref{2.3}) is
bounded by
\begin{equation}
\label{2.4} {(2\mu)^{2N} \over (N+2)!} X \mu + {A^{B^N}
(2\mu)^{2N} \over (2N+3)!} {X \mu \over \log{X}}.
\end{equation}
Using (\ref{1.6}), the main term of (\ref{2.3}) contributes
\begin{equation}
\label{2.5} {X \over 2G_{\lambda}} \int_{0}^{\mu} {1 \over
u^{1+\lambda}} \sum_{j=1}^{N} {(-1)^{j+1} \over (2j)!}
{\Gamma(2j+1) \over \Gamma(j+1)} (2u)^{2j} du + O\Bigl({X \over
\log{X}} \sum_{j=1}^{N} {A^{B^j} \over (2j)!} \int_{0}^{\mu}
u^{2j-2} du \Bigr).
\end{equation}
The above error term is bounded by
\begin{equation}
\label{2.6} {X \over \log{X}} A^{B^N} \sum_{j=1}^{N} {\mu^{2j-1}
\over (2j)!} \ll {X \over \log{X}} A^{B^N} e^{\mu}.
\end{equation}
Using (\ref{2.1}) with $z=j+{1 \over 2}$, the main term of
(\ref{2.5}) becomes
\begin{eqnarray*}
& &{X \over 2G_{\lambda}} \int_{0}^{\mu} {1 \over u^{1+\lambda}}
\sum_{j=1}^{N} {(-1)^{j+1} \over (2j)!} {2^{2j} \over \sqrt{\pi}}
\Gamma(j+{1 \over 2})
(2u)^{2j} du \\
&=&{X \over 2 \sqrt{\pi} G_{\lambda}} \int_{0}^{\mu} {1 \over
u^{1+\lambda}} \sum_{j=1}^{N} {(-1)^{j+1} \over (2j)!} (4u)^{2j}
\Bigl(\int_{0}^{\infty}
x^{j-{1/2}} e^{-x} dx\Bigr) du \\
&=&{X \over \sqrt{\pi} G_{\lambda}} \int_{0}^{\infty} x^{-{1/2}}
e^{-x} \int_{0}^{\mu} {1 \over u^{1+\lambda}}
\Bigl((\sin{2\sqrt{x} u})^2 + O\bigl( {(4\sqrt{x} u)^{2N+2} \over
(2N+2)!}\bigr) \Bigr)du dx.
\end{eqnarray*}
The contribution from the above error term is
\begin{equation}
\label{2.7} \ll_{\lambda} {4^{2N+2} (N+1)! \over (2N+2)!}
\mu^{2N+2-\lambda} X.
\end{equation}
The contribution from the main term is
\begin{eqnarray*}
&=&{X \over \sqrt{\pi} G_{\lambda}} \int_{0}^{\infty} x^{-{1/2}}
e^{-x} \Bigl(\int_{0}^{\infty} {(\sin{2\sqrt{x} u})^2 \over
u^{1+\lambda}} du + O({1 \over \lambda \mu^{\lambda}}) \Bigr) dx \\
&=&{2^{\lambda} \over \sqrt{\pi}} \Gamma({\lambda + 1 \over 2}) X
+ O_{\lambda} \Bigl({X \over \mu^{\lambda}} \Bigr) \\
&=& {\Gamma(\lambda+1) \over \Gamma({\lambda \over 2}+1)} X +
O_{\lambda} \Bigl({X \over \mu^{\lambda}} \Bigr)
\end{eqnarray*}
by the definition of $G_{\lambda}$ and (\ref{2.1}). Combining this
with (\ref{2.2}), (\ref{2.3}), (\ref{2.4}), (\ref{2.5}),
(\ref{2.6}) and (\ref{2.7}), we have
\begin{equation}
\label{2.8} \int_{1}^{X} |F(x)|^{\lambda} dx = {\Gamma(\lambda+1)
\over \Gamma({\lambda \over 2}+1)} X + \mbox{error},
\end{equation}
where
$$\mbox{error} \ll_{\lambda} {X \over \mu^{\lambda}} + {(2\mu)^{2N+1} \over
(N+1)!} X + {A^{B^N} (2\mu)^{2N+2} \over (2N+3)!} {X \over
\log{X}} + A^{B^N} e^{\mu} {X \over \log{X}}.$$ Now, we choose $N
= {\log_{3}{X} \over \log_{4}{X}}$ and $\mu = \sqrt[\lambda]
{\log_{4}{X}}$, then one can check that $A^{B^N} \ll
\sqrt{\log{X}}$, $e^{\mu} \ll_{\lambda} (\log{X})^{1 \over 4}$,
and using Stirling's formula,
\begin{eqnarray*}
{(2\mu)^{2N+2} \over (N+1)!} &\ll& {e^{(2N+2)\log{2\mu}} \over
\sqrt{N} (N/e)^
{N}} \ll {1 \over \sqrt{N}} e^{(2N+2)\log{2\mu} - N\log{(N/e)}} \\
&\ll_{\lambda}& {1 \over \sqrt{N}} \ll {1 \over \log_{3}{X}}.
\end{eqnarray*}
Consequently, (\ref{2.8}) becomes
$$\int_{1}^{X} |F(x)|^{\lambda} dx = {\Gamma(\lambda+1) \over \Gamma({\lambda
\over 2}+1)} X + O_{\lambda}\Bigl({X \over \log_{4}{X}}\Bigr)
$$
which gives the theorem for $0 < \lambda \leq 1$ after multiplying
through by $({1 \over \sqrt{2}} h^{1/2} (\log{X \over
h})^{1/2})^{\lambda}$.
%%-----------------------------------------------------------------------------

(ii) $1 < \lambda$.

Let $\lambda = 2m+1+\theta$ where $m$ is a non-negative integer
and $0 < \theta \leq 2$. Let $D_{\theta} = \int_{0}^{\infty}
{(\sin{u})^4 \over u^{2+\theta}} du$. Since
\begin{eqnarray*}
|F(x)|^{\lambda} &=& {|F(x)|^{2m} \over D_{\theta}}
\int_{0}^{\infty}
{(\sin{|F(x)|u})^4 \over u^{2+\theta}} du \\
&=& {|F(x)|^{2m} \over D_{\theta}} \int_{0}^{Y} {(\sin{|F(x)|u})^4
\over u^{2+\theta}} du + O_{\lambda}(|F(x)|^{2m} Y^{-1-\theta}).
\end{eqnarray*}
Then, similar to the calculation in case (i), for some $Y \geq 1$,
\begin{equation}
\label{2.9} \int_{1}^{X} |F(x)|^{\lambda} dx = {1 \over
D_{\theta}} \int_{0}^{Y} {1 \over u^{2+\theta}} \Bigl(\int_{1}^{X}
|F(x)|^{2m} (\sin{|F(x)|u})^4 dx \Bigr) du +
O_{\lambda}(Y^{-1-\theta} X)
\end{equation}
by (\ref{1.6}). From ($7$) of Ghosh [\ref{G2}], we have
$$(\sin{u})^4 = {1 \over 8}\sum_{j=2}^{N} b_j u^{2j} + O\Bigl({(4u)^{2N+2}
\over (2N+2)!} \Bigr) \mbox{ where } b_j = {(-4)^{j+1} \over
(2j)!} (4^{j-1} - 1),$$ and $N$ is an integer, exceeding $2$,
which will be chosen later. Using this Taylor series and
(\ref{1.6}), the main term of (\ref{2.9}) equals
\begin{equation}
\label{2.10} {X \over 8D_{\theta}} \int_{0}^{Y} {1 \over
u^{2+\theta}} \sum_{j=2}^{N} b_j u^{2j} {\Gamma(2(m+j)+1) \over
\Gamma((m+j)+1)} du
\end{equation}
$$+O\Bigl({X \over \log{X}} \sum_{j=2}^{N} A^{B^j} |b_j| \int_{0}^{Y} u^{2j-2-
\theta}du\Bigr) + \int_{1}^{X} |F(x)|^{2m} {(4|F(x)|)^{2N+2} \over
(2N+2)!} \int_{0}^{Y} u^{2N-\theta} du dx \Bigr).$$ The error term
of (\ref{2.10}) is
\begin{equation}
\label{2.11} \ll {A^{B^N} e^{4Y} \over Y^{1+\theta}}{X \over
\log{X}} + {4^{2N+2} Y^{2N + 1} \over (2N+2)!} \Bigl({(2(m+N+1))!
\over (m+N+2)!} + {A^{B^{m+N+1}} \over \log{X}} \Bigr) X.
\end{equation}
By (\ref{2.1}), the main term of (\ref{2.10})
\begin{eqnarray*}
&=& {X \over D_{\theta}} \int_{0}^{Y} {1 \over u^{2+\theta}}
\sum_{j=2}^{N} {b_j \over 8} u^{2j} {2^{2(m+j)} \over \sqrt{\pi}}
\Bigl(\int_{0}^{\infty}
x^{m+j-{1/2}} e^{-x} dx \Bigr) du \\
&=& {2^{2m} X \over \sqrt{\pi} D_{\theta}} \int_{0}^{\infty}
x^{m-{1/2}} e^{-x} \int_{0}^{Y} {1 \over u^{2+\theta}}
\Bigl((\sin{2\sqrt{x} u})^4 + O({ (8\sqrt{x} u)^{2N+2} \over
(2N+2)!}) \Bigr) du dx.
\end{eqnarray*}
Contribution from the error is
\begin{equation}
\label{2.12} \ll_{\lambda} {(N+m+1)! \over (2N+2)!} 2^{4N}
Y^{2N+1} X
\end{equation}
while the main term
\begin{eqnarray*}
&=&{2^{2m} X \over \sqrt{\pi} D_{\theta}} \int_{0}^{\infty}
x^{m-{1/2}} e^{-x} \Bigl(\int_{0}^{\infty} {(\sin{2\sqrt{x}u})^4
\over u^{2+\theta}} du +
O(Y^{-1-\theta}) \Bigr) dx \\
&=&{2^{2m+1+\theta} \over \sqrt{\pi}} X \int_{0}^{\infty}
x^{m+{\theta/2}
} e^{-x} dx + O_{\lambda}(X Y^{-1-\theta}) \\
&=&{\Gamma(\lambda+1) \over \Gamma({\lambda \over 2}+1)} X +
O_{\lambda}(X Y^ {-1-\theta})
\end{eqnarray*}
by (\ref{2.1}). Therefore, combining this with (\ref{2.9}),
(\ref{2.10}), (\ref{2.11}) and (\ref{2.12}), we have
\begin{equation}
\label{2.13} \int_{1}^{X} |F(x)|^{\lambda} dx = {\Gamma(\lambda+1)
\over \Gamma({\lambda \over 2}+1)} X + \mbox{error},
\end{equation}
where
\begin{equation}
\label{2.14}
\begin{split}
\mbox{error} &\ll_{\lambda} X \Bigl({1 \over Y} + {A^{B^N} e^{4Y}
\over Y
\log{X}} +4^{2N} Y^{2N+1} {(2(m+N+1))! \over (2N+2)! (m+N+2)!} \\
&+ {4^{2N} A^{B^N} \over (2N+2)!}{Y^{2N+1} \over \log{X}} +
{4^{2N} (N+m+1)! \over (2N+2)!} Y^{2N+1} \Bigr).
\end{split}
\end{equation}
By Stirling's formula, (\ref{2.14})
$$\ll_{\lambda} X \Bigl({1 \over Y} + {A^{B^N} e^{4Y} \over Y \log{X}} +
{4^{4N} \over (N+1)^{N-m}} Y^{2N+1} + {A^{B^N} \over (N+1)^{2N+2}}
{Y^{2N+1} \over \log{X}} \Bigr).$$ Now, pick $N = {\log_{3}{X}
\over \log_{4}{X}}$ and $Y = \sqrt{\log_{3}{X}}$, we have $A^{B^N}
\ll \sqrt{\log{X}}$, $4^{4N} \ll \sqrt[4]{\log{\log{X}}}$,
$(N+1)^{N-m} \gg_{\lambda} \log{\log{X}}$, $e^{4Y} \ll
\log{\log{X}}$ and $Y^{2N+1} \ll \sqrt{\log{\log{X}}}$. Thus,
(\ref{2.14}) is $O_{\lambda}({X \over \sqrt[3]{\log_{3}{X}}})$ and
we get the theorem after multiplying (\ref{2.13}) by $({1 \over
\sqrt{2}} h^ {1/2} (\log{X \over h})^{1/2})^{\lambda}$.
%%----------------------------------------------------------------------------
\section{``Equivalence'' for odd moments}

In [\ref{C}], the author proved that (\ref{1.1}) and (\ref{1.3})
are roughly equivalent to one another when $k$ is even. One would
like to prove a similar statement when $k$ is odd. However, the
difficulty lies in that we no longer have asymptotic formulas.
But, if one has good knowledge about all the even moments then it
is possible to get the following
\begin{thm}
\label{theorem3} Assume Riemann Hypothesis. If (\ref{1.6}) holds 
in $X^{\epsilon} \leq h \leq X^{1-\epsilon}$ for some $\epsilon > 0$
and all positive even integer $k$, then, for any positive odd integer $n$,
$$\int_{1}^{X} (\psi(x+h)-\psi(x)-h)^n dx = o(X h^{n/2} (\log{X})^{n/2})$$
for $X^{\epsilon} \leq h \leq X^{1-\epsilon}$ implies that, for
some $\epsilon_1 > 0$,
$$\int_{1}^{X} (\psi(x+\delta x)-\psi(x)-\delta x)^n dx = o\Bigl(X^{{n/2}+1}
\delta^{n/2} (\log{1 \over \delta})^{n/2}\Bigr)$$ for $X^{-1
+2\epsilon + 2\epsilon_1} \leq \delta \leq X^{-\epsilon}/2$.
\end{thm}

Conversely, one also has
\begin{thm}
\label{theorem4} Assume Riemann Hypothesis. If (\ref{1.7}) holds 
in $X^{-1+\epsilon} \leq \delta \leq X^{-\epsilon}$ for some 
$\epsilon > 0$ and all positive even integer $k$, then, for any
positive odd integer $n$,
$$\int_{1}^{X} (\psi(x+\delta x)-\psi(x)-\delta x)^n dx = o\Bigl(X^{{n/2}+1}
\delta^{n/2} (\log{1 \over \delta})^{n/2}\Bigr)$$ for
$X^{-1+\epsilon} \leq \delta \leq X^{-\epsilon}$ implies that, for
some $\epsilon_1 > 0$,
$$\int_{1}^{X} (\psi(x+h)-\psi(x)-h)^n dx = o(X h^{n/2} (\log{X})^{n/2})$$
for $X^{2\epsilon + \epsilon_1} \leq h \leq X^{1-(n/2+1) \epsilon
- 2\epsilon_1} / 2$.
\end{thm}

Remarks: 1. The proofs of the above theorems are very similar to
the proofs of theorems in [\ref{C}]. We shall give a sketch for
Theorem \ref{theorem3} only. 2. We did not optimize the ranges for
$h$ and $\delta$ in the above theorems. Improvements are possible
since the error terms in the proofs in [\ref{C}] are smaller.

\bigskip

Sketch of proof of Theorem \ref{theorem3}: Observe that
\begin{eqnarray*}
& &\int_{1}^{X} |\psi(x+h)-\psi(x)-h|^n dx + \int_{1}^{X}
(\psi(x+h)-\psi(x)-h)
^n dx \\
&=& 2 \int_{1}^{X} \max \{\psi(x+h)-\psi(x)-h, 0\}^n dx.
\end{eqnarray*}
Thus, by the assumptions in Theorem \ref{theorem3}, we can apply
Theorem \ref{theorem1} and get
\begin{equation}
\label{3.1}
\begin{split}
\int_{1}^{X} \max \{\psi(x+h)-\psi(x)-h, 0\}^n dx &= {1 \over 2}
{\Gamma(n+1) \over \Gamma({n \over 2}+1) 2^{n/2}} X h^{n/2}
(\log{X \over h})^{n/2} \\
&+ o(X h^{n/2} (\log{X})^{n/2})
\end{split}
\end{equation}
Now, for any $0 < \epsilon_1 < \epsilon$, one can imitate Saffari
and Vaughan's argument as in Theorem 3.1 of [\ref{C}] and get
\begin{eqnarray*}
& &\int_{0}^{\Delta} \int_{1}^{X} max\{\psi(x+\delta
x)-\psi(x)-\delta x,0\}^n dx \; d\delta \\
&=&{\Gamma(n+1) \over 2 \Gamma({n \over 2}+2) 2^{n/2}}
\int_{0}^{\Delta X} h^{n/2} \Bigl(\log{X \over h}\Bigr)^{n/2} dh +
o(\Delta^{{n/2}+1} X^{{n/2}+1} (\log{1 \over \Delta})^{n/2})
\end{eqnarray*}
for $X^{-1 +2\epsilon + \epsilon_1} \leq \delta \leq
X^{-\epsilon}$. Let $f(x,h)=max\{\psi(x+h)-\psi(x)-h, 0\}$,
$g(x,\delta x)=f(x,\Delta x)$ for $\Delta \leq \delta \leq
(1+\eta) \Delta$. Following the argument in Theorem 3.1 of
[\ref{C}] without choosing $\eta$ explicitly or simply following
the argument in [\ref{GM}], one has, for $X^{-1+ 2\epsilon +
2\epsilon_1} \leq \Delta \leq X^{-\epsilon}/2$,
\begin{equation}
\label{3.2}
\begin{split}
& \int_{1}^{X} \max \{\psi(x+\Delta x) - \psi(x) - \Delta x,0\}^n
dx \\
=& {1 \over 2} {\Gamma(n+1) \over \Gamma({n \over 2}+2) 2^{n/2}}
X^{{n/2}+1} \Delta^{n/2} (\log{1 \over \Delta})^{n/2} +
o\Bigl(X^{{n/2}+1} \Delta^{n/2} (\log{1 \over \Delta})^{n/2}\Bigr)
\end{split}
\end{equation}
by letting $\eta$ approach $0$ sufficiently slowly. The only
difference in the argument is the use of
\begin{eqnarray*}
|f(x,\delta x)-g(x,\delta x)| &\leq& |(\psi(x+\delta
x)-\psi(x)-\delta x) -
(\psi(x+\Delta x)-\psi(x)-\Delta x)| \\
&=&|\psi\bigl((x+\Delta x)+(\delta-\Delta)x\bigr) -
(\delta-\Delta)x|
\end{eqnarray*}
for $\Delta \leq \delta \leq (1+\eta) \Delta$ when one estimates
the integral $\int \int |f-g|^n$. The above is justified by
$$|\max\{a,0\} - \max\{b,0\}| \leq |a-b|$$
which can be easily verified by considering different cases of
signs for $a$ and $b$. Similarly, one has
\begin{equation}
\label{3.3}
\begin{split}
& \int_{1}^{X} \min\{\psi(x+\Delta x) - \psi(x) - \Delta x,0\}^n
dx \\
=& -{1 \over 2}{\Gamma(n+1) \over \Gamma({n \over 2}+2) 2^{n/2}}
X^{{n/2}+1} \Delta^{n/2} (\log{1 \over \Delta})^{n/2} +
o\Bigl(X^{{n/2}+1} \Delta^{n/2} (\log{1 \over
\Delta})^{n/2}\Bigr).
\end{split}
\end{equation}
Consequently, adding (\ref{3.2}) and (\ref{3.3}), we have the
theorem.
%%---------------------------------------------------------------------------
\section{Numerical evidence}

Instead of having the first main terms only in Theorem
\ref{theorem1} and \ref{theorem2}, one should expect more to be
true, namely,
\begin{conj}
\label{conj1} For every $\epsilon > 0$ and $\lambda > 0$,
$$\int_{1}^{X} |\psi(x+h)-\psi(x)-h|^{\lambda} dx = {\Gamma(\lambda+1) \over
\Gamma({\lambda \over 2}+1) 2^{\lambda/2}} h^{\lambda/2+1}
\int_{E}^{X/h} \bigl(\log{x \over E}\bigr)^{\lambda/2} dx + o(X
h^{\lambda/2})$$ for $X^{\epsilon} \leq h \leq X^{1-\epsilon}$.
Here $E = 2\pi e^{C_0 -1}$.
\end{conj}
\begin{conj}
\label{conj2} For every $\epsilon > 0$ and $\lambda > 0$,
$$\int_{1}^{X} |\psi(x+\delta x) - \psi(x) - \delta x|^{\lambda} dx = {\Gamma
(\lambda+1) \over \Gamma({\lambda \over 2} + 2) 2^{\lambda/2}}
X^{{\lambda/2}+1 } \delta^{\lambda/2} \Bigl(\log{1 \over
E\delta}\Bigr)^{\lambda/2} + o(X^{{\lambda/2}+1}
\delta^{\lambda/2})$$ for $X^{-1+\epsilon} \leq \delta \leq
X^{-\epsilon}$. Again $E = 2\pi e^{C_0 -1}$.
\end{conj}
Using a C program, we get some numerical evidence in support of
Conjecture \ref{conj2} as well as the odd moments for (\ref{1.2}).
Firstly, regarding Conjecture \ref{conj2}:

\bigskip
For $X=10^8$ and $\delta=10^{-4}$:
\bigskip

\begin{math}
\begin{array}{lll}
\lambda & \mbox{Actual value for $\lambda$-th moment} &
\mbox{Result from formula} \\
1.0 & 1.5009\ast 10^{10} & 1.4851\ast 10^{10} \\
2.1 & 7.1441\ast 10^{12} & 6.9344\ast 10^{12} \\
3.2 & 4.8737\ast 10^{15} & 4.6213\ast 10^{15} \\
4.3 & 4.1913\ast 10^{18} & 3.8864\ast 10^{18} \\
5.4 & 4.2519\ast 10^{21} & 3.8768\ast 10^{21} \\
6.5 & 4.8884\ast 10^{24} & 4.4213\ast 10^{24}
\end{array}
\end{math}

\bigskip
For $X=10^{10}$ and $\delta=10^{-5}$:
\bigskip

\begin{math}
\begin{array}{lll}
\lambda & \mbox{Actual value for $\lambda$-th moment} &
\mbox{Result from formula} \\
1.0 & 5.3464\ast 10^{12} & 5.3452\ast 10^{12} \\
2.1 & 1.0218\ast 10^{16} & 1.0210\ast 10^{16} \\
3.2 & 2.7871\ast 10^{19} & 2.7835\ast 10^{19} \\
4.3 & 9.5892\ast 10^{22} & 9.5764\ast 10^{22} \\
5.4 & 3.9120\ast 10^{26} & 3.9079\ast 10^{26} \\
6.5 & 1.8248\ast 10^{30} & 1.8232\ast 10^{30}
\end{array}
\end{math}

\bigskip
Note: We just happen to pick some values for $\lambda$.
\bigskip

Secondly, for $n$ odd,
$$\int_{1}^{X} \bigl(\psi(x+\delta x) - \psi(x) - \delta x \bigr)^n dx =
o\Bigl(X^{{n/2}+1} \delta^{n/2} (\log{1 \over
\delta})^{n/2}\Bigr).$$

\bigskip
For $X=10^8$ and $\delta=10^{-4}$:
\bigskip

\begin{math}
\begin{array}{lll}
n & \mbox{Actual value for $n$-th moment} & {\Gamma(n+1) \over
\Gamma({n \over 2}+2) 2^{n/2}} X^{{n/2}+1} \delta^{n/2}
(\log{1 \over \delta})^{n/2} \\
1 & -4.9574\ast 10^{7} & 1.6143\ast 10^{10} \\
3 & -2.0632\ast 10^{13} & 1.7842\ast 10^{15} \\
5 & -3.3174\ast 10^{18} & 4.6952\ast 10^{20}
\end{array}
\end{math}

\bigskip
For $X=10^{10}$ and $\delta=10^{-5}$:
\bigskip

\begin{math}
\begin{array}{lll}
n & \mbox{Actual value for $n$-th moment} & {\Gamma(n+1) \over
\Gamma({n \over 2}+2) 2^{n/2}} X^{{n/2}+1} \delta^{n/2}
(\log{1 \over \delta})^{n/2} \\
1 & 7.2371\ast 10^{8} & 5.7074\ast 10^{12} \\
3 & -1.3468\ast 10^{16} & 7.8851\ast 10^{18} \\
5 & -2.5587\ast 10^{23} & 2.5937\ast 10^{25}
\end{array}
\end{math}

\bigskip
Note: ${\Gamma(n+1) \over \Gamma({n/2}+2) 2^{n/2}}$ acts as a
normalization constant coming from the main term of $\lambda$-th
moment.
%-----------------------------------------------------------------

%%-----------------------------------------------------------------------------
Tsz Ho Chan\\
American Institute of Mathematics\\
360 Portage Avenue\\
Palo Alto, CA 94306\\
USA\\
thchan@aimath.org

\end{document}